\documentclass[12pt,reqno]{amsart}
\usepackage{amsmath,amssymb,amsthm}
\usepackage[margin=1in]{geometry}
\usepackage[hidelinks]{hyperref}
\newtheorem{theorem}{Theorem}
\newtheorem{corollary}{Corollary}
\newtheorem{lemma}{Lemma}
\theoremstyle{remark}
\newtheorem{remark}{Remark}
\newcommand{\eps}{\varepsilon}

\newcommand{\Tam}{\operatorname{Tam}}
\DeclareFontFamily{U}{wncy}{}
\DeclareFontShape{U}{wncy}{m}{n}{<->wncyr10}{}
\DeclareSymbolFont{mcy}{U}{wncy}{m}{n}
\DeclareMathSymbol{\Sha}{\mathord}{mcy}{"58}
\begin{document}
\title[Extreme $L$-values of almost prime quadratic twists II]{Extreme central $L$-values of almost prime quadratic twists of elliptic curves II}
\author{Shenghao Hua}
\address{Shanghai Institute for Mathematics and Interdisciplinary Sciences (SIMIS), Shanghai 200433, China}
\address{Research Institute of Intelligent Complex Systems, Fudan University, Shanghai, 200433, China}
\email{huashenghao@vip.qq.com}

\date{}
\begin{abstract}
We prove lower bounds for central $L$-values of quadratic twists of an elliptic curve by discriminants with at most five prime factors. Under the rank-zero Birch--Swinnerton-Dyer formula, these bounds give large Tate--Shafarevich groups.
\end{abstract}
\maketitle
\section{Introduction}

The resonance method finds extreme values of $L$-functions by weighting them with suitable Dirichlet polynomials. Soundararajan \cite{Sound} used it to study large values of the Riemann zeta-function and extreme central values in families of $L$-functions. We use the method to find large values in a sparse family of quadratic twists.

For quadratic twists of elliptic curves, the resonance method gives large central $L$-values \cite{HL}.
Hua and Huang \cite{HH} proved such bounds for discriminants with at most $W$ prime factors, for every fixed $W\ge20$.
Here we obtain bounds for $W\ge5$.

Let $E/\mathbb Q$ be an elliptic curve of conductor $N$, and let
$E_d$ be its quadratic twist by $d$. We normalize the $L$-function
so that its central point is $s=1/2$.
Put $N_0=\operatorname{lcm}(8,N)$.
A class $a\pmod{N_0}$ is admissible if $(a,N_0)=1$ and
$a\equiv1$ or $5\pmod8$. For such a class, let $\mathcal D_a^-$ be the set of negative fundamental discriminants $d\equiv a\pmod{N_0}$.
Write $w_a$ for the common root number of these twists.
We fix an admissible class with $w_a=+1$.
For an integer $W\ge1$, let $\mathcal F_{a,W}(X)$ be the set of $d\in\mathcal D_a^-$ with $X\le|d|\le2X$ and $\omega(|d|)\le W$, where $\omega(n)$ counts the distinct prime divisors of $n$. Put $H_X=\sqrt{\log X/\log\log X}$.

\begin{theorem}\label{thm:main}
Let $E$ and $a\pmod{N_0}$ be as above.
For every fixed integer $W\ge5$ and every fixed
$0<c<\sqrt{(7W-34)/(7W+1)}$, the set $\mathcal F_{a,W}(X)$
is nonempty for all sufficiently large $X$, and
\[
 \max_{d\in\mathcal F_{a,W}(X)}L(1/2,E_d)\ge\exp(cH_X).
\]
\end{theorem}

For $W=5$, one may take any $0<c<1/6$.
The proof combines the resonance method with a weighted sieve.

Large central values also give large Tate--Shafarevich groups under the Birch--Swinnerton-Dyer formula. Goldfeld and Szpiro \cite{GS} studied bounds for their orders, and Radziwi\l\l\ and Soundararajan \cite{RS} studied their distribution in quadratic twist families. When $L(1/2,E_d)\ne0$, the group $\Sha(E_d/\mathbb Q)$ is finite \cite{GZ,Kol}. The rank-zero formula relates its order to the central value, the real period, the rational torsion and the Tamagawa product. For the twists here, the period is proportional to $|d|^{-1/2}$, and the Tamagawa product is bounded when $\omega(|d|)\le W$. Hua and Huang \cite{HH} used these facts to find large Tate--Shafarevich groups in almost-prime twist families. Theorem~\ref{thm:main} gives the following result with at most five prime factors.

\begin{corollary}\label{cor:Sha-main}
Let $E$, $a$, $W$ and $c$ be as in Theorem~\ref{thm:main}.
Assume the rank-zero Birch--Swinnerton-Dyer formula for the twists
$E_d$ with $d\in\mathcal D_a^-$ and $L(1/2,E_d)\ne0$.
For every sufficiently large $X$, there is
$d\in\mathcal F_{a,W}(X)$ such that
\begin{equation}\label{eq:Sha-large}
 \#\Sha(E_d/\mathbb Q)\ge\sqrt X\exp(cH_X).
\end{equation}
For elliptic curves with complex multiplication the conclusion holds
unconditionally.
\end{corollary}
\section{The first moment and the resonator}\label{sec:resonator}

Write $L(s,E)=\sum_{n\ge1}a(n)n^{-s}$ and
$\chi_d(n)=(\frac dn)$.
The central values $L_d=L(1/2,E_d)$ are nonnegative \cite{RS}.
Fix a nonzero nonnegative $\Phi\in C_c^\infty((1,2))$.

For $p\nmid N_0$, write
$Q_p=1-(a(p)^2-2)/p+p^{-2}$ and $T_p=1+p^{-1}+Q_p/p$.
Define $g$ multiplicatively by $g(1)=1$, $g(p^{2r-1})=T_p^{-1}$ and $g(p^{2r})=(1+p^{-1})T_p^{-1}$ for $r\ge1$.
Put $\delta(p)=1-g(p^2)=Q_p/(pT_p)$, and extend $\delta$
multiplicatively to squarefree integers coprime to $N_0$.
By \cite[Lemma~2.1, $r=0$]{HL}, for $(t,N_0)=1$ and
$t=t_1t_2^2$ with $t_1$ squarefree,
\[
 \sum_{d\in\mathcal D_a^-}L_d\chi_d(t)\Phi(|d|/X)
 =C_1X\frac{a(t_1)}{\sqrt{t_1}}g(t)
   +O_{E,\Phi,\eps}(X^{1/2+\eps}t^{1/2+\eps}),
\]
where $C_1>0$ depends only on the fixed data.
Write
\[
 \mathcal S_a(0;X,u,\ell)
 =\sum_{\substack{d\in\mathcal D_a^-\\\ell\mid d}}
                    L_d\chi_d(u)\Phi(|d|/X).
\]

Fix $\theta>0$, $0<\rho<1$, and put
$M=X^\theta$, $\mathcal L=\rho\sqrt{\log M\log\log M}$,
$P_- =\mathcal L^2$, and $P_+=\exp((\log\mathcal L)^2)$.
Define a real multiplicative function $b$ supported on squarefree integers by $b(1)=1$. For primes $P_-\le p\le P_+$ with $p\nmid N_0$, set $b(p)=a(p)\mathcal L/(\sqrt p\log p)$; set $b(p)=0$ for all other primes.
Let $R(d)=\sum_{m\le M}b(m)\chi_d(m)$, $A_0(p)=1+b(p)^2p/(p+1)$ and $A_1(p)=1+b(p)^2g(p^2)+2a(p)b(p)g(p)/\sqrt p$.
Let $B_i=\prod_{p\nmid N_0}A_i(p)$ for $i=0,1$.
These products are finite. Rankin--Selberg theory and partial
summation give
\begin{equation}\label{eq:sv-prime-sums}
 \sum_{P_-\le p\le P_+}\frac{a(p)^2}{p\log p}
 =\frac{1+o(1)}{2\log\mathcal L},\qquad
 \sum_{P_-\le p\le P_+}\frac{a(p)^2}{p}
 =O(\log\log\mathcal L);
\end{equation}
see \cite[Lemma~3]{RS}. Since $P_-\to\infty$, the primes in this interval do not divide $N_0$ for large $X$.

\begin{lemma}
For the coefficients above,
\begin{equation}\label{eq:sv-coefficient-bound}
 \sum_{m\le M}|b(m)|\sqrt m\le M^{1+o(1)}.
\end{equation}
There is $c_\rho>0$ such that, with
$\tau_M=\exp(-c_\rho\log M/(\log\log M)^3)$, uniformly for
squarefree $\ell$ coprime to $N_0$,
\begin{equation}\label{eq:sv-truncation-first}
 \sum_{\substack{m_1,m_2\le M\\ (m_1m_2,\ell)=1}}
 b(m_1)b(m_2)g(m_1m_2)
 \frac{a((m_1m_2)_1)}{\sqrt{(m_1m_2)_1}}
 =\prod_{p\nmid\ell N_0}A_1(p)(1+O(\tau_M)).
\end{equation}
Here $(m_1m_2)_1$ is the squarefree part of $m_1m_2$.
\end{lemma}
\begin{proof}
Hasse's bound $|a(p)|\le2$ and the elementary harmonic estimate give
\[
 \sum_{m\le M}|b(m)|\sqrt m
 \le M\prod_{P_-\le p\le P_+}
             \left(1+\frac{|a(p)|\mathcal L}{p\log p}\right)
 \le M\exp\!\left(4\mathcal L(1+(\log\mathcal L)^2)\right)
 =M^{1+o(1)}.
\]
Here $\log p\ge\log2$, and the harmonic sum bounds the exponent sum by $4(1+\log P_+)$ for large $M$. Also $\mathcal L(1+(\log\mathcal L)^2)=o(\log M)$.
For the second assertion write $m_1=qr$, $m_2=qs$, with
$q,r,s$ pairwise coprime and squarefree. The summand becomes
\[
 b(q)^2g(q^2)
 \frac{a(r)b(r)g(r)}{\sqrt r}
 \frac{a(s)b(s)g(s)}{\sqrt s},
\]
which is nonnegative. The unrestricted sum is the product of
$A_1(p)$ over $p\nmid\ell N_0$.
Put $\nu=(\log\mathcal L)^{-3}$.
Rankin's trick bounds the relative contribution of $qr>M$ or
$qs>M$ by
\[
 2\exp\left(-\nu\log M+
 \sum_{p\nmid\ell N_0}
 \left(b(p)^2g(p^2)+\frac{a(p)b(p)g(p)}{\sqrt p}\right)
 (p^\nu-1)\right).
\]
On the support, $p^\nu-1=\nu\log p(1+o(1))$ uniformly.
By \eqref{eq:sv-prime-sums}, the quadratic term with the factor
$\log p$ sums to $(\rho^2+o(1))\log M$.
The linear term sums to
$O(\mathcal L\log\log\mathcal L)=o(\log M)$.
Omitting primes dividing $\ell$ decreases these sums.
The exponent is at most $-(1-\rho^2+o(1))\nu\log M$.
\end{proof}

\begin{lemma}\label{lem:euler-ratio}
We have $\log(B_1/B_0)=(2\rho\sqrt\theta+o(1))H_X$.
\end{lemma}
\begin{proof}
On the support of $b$, put
$x_p=2a(p)b(p)g(p)/\sqrt p$ and $y_p=b(p)^2g(p^2)$.
Then
\[
 \log\frac{A_1(p)}{A_0(p)}
 =\log\frac{1+y_p}{1+b(p)^2p/(p+1)}
  +\log\left(1+\frac{x_p}{1+y_p}\right).
\]
Since $g(p),g(p^2)=1+O_E(p^{-1})$ and
$b(p)^2\ll(\log\mathcal L)^{-2}$, this is
\[
 \frac{2a(p)b(p)g(p)}{\sqrt p}
 +O\left(\frac{b(p)^2}{p}+\frac{a(p)^2b(p)^2}{p}
              +\frac{|a(p)||b(p)|^3}{\sqrt p}\right).
\]
Partial summation gives total error
$O(\mathcal L/(\log\mathcal L)^4)=o(H_X)$.
Replacing $g(p)$ by $1$ in the linear term contributes $o(1)$.
By \eqref{eq:sv-prime-sums},
\[
 \log\frac{B_1}{B_0}
 =(1+o(1))\frac{\mathcal L}{\log\mathcal L}
 =(2\rho\sqrt\theta+o(1))H_X.
\]
\end{proof}

\section{Divisor sums}
\begin{lemma}\label{lem:richert-weak-divisor}
Let $\ell$ be squarefree and $(u,\ell)=(u\ell,N_0)=1$.
Writing $u=u_1u_2^2$ with $u_1$ squarefree, we have
\[
 \mathcal S_a(0;X,u,\ell)
 =C_1X\frac{a(u_1)}{\sqrt{u_1}}g(u)\delta(\ell)
 +O\bigl(X^{1/2+\eps}u^{1/2+\eps}\ell^{1+\eps}\bigr).
\]
\end{lemma}
\begin{proof}
For a fundamental discriminant $d$,
$\mathbf1_{\ell\mid d}=\sum_{q\mid\ell}\mu(q)\chi_d(q^2)$.
Apply \cite[Lemma~2.1, $r=0$]{HL} to each insertion $uq^2$.
Since $(u,\ell)=1$, all squarefree parts are $u_1$ and $g(uq^2)=g(u)g(q^2)$. As $\ell$ is squarefree, $\sum_{q\mid\ell}\mu(q)g(q^2)=\delta(\ell)$. Thus $\sum_{q\mid\ell}\mu(q)g(uq^2)=g(u)\delta(\ell)$.
The errors sum to
$O(X^{1/2+\eps}u^{1/2+\eps}\sum_{q\mid\ell}q^{1+2\eps})$.
Use $\sum_{q\mid\ell}q^{1+2\eps}\ll_\eps\ell^{1+3\eps}$
and choose smaller intermediate exponents.
\end{proof}

\begin{lemma}\label{lem:richert-distribution}
Let $M=X^\theta$, $D=X^\gamma$ and $4\theta+4\gamma<1$.
Put $\mathcal X=C_1XB_1$, $G(p)=\delta(p)/A_1(p)$ for
$p\nmid N_0$, and $G(p)=0$ for $p\mid N_0$.
Extend $G$ multiplicatively to squarefree integers.
For squarefree $\ell$ coprime to $N_0$, write
\[
 \mathcal A_\ell
 =\sum_{\substack{d\in\mathcal D_a^-\\\ell\mid d}}
 L(1/2,E_d)R(d)^2\Phi(|d|/X)
 =\mathcal XG(\ell)+r_\ell.
\]
Set $r_\ell=0$ for all other $\ell$.
Then, for every fixed $A>0$,
\[
 \sum_{\substack{\ell\le D\\\mu(\ell)^2=1\\ (\ell,N_0)=1}}
 |r_\ell|\ll_A\mathcal X(\log X)^{-A}.
\]
\end{lemma}
\begin{proof}
Expand $R(d)^2$. Terms with $(m_1m_2,\ell)>1$ vanish.
Apply Lemma~\ref{lem:richert-weak-divisor} to the other terms.
The truncation formula \eqref{eq:sv-truncation-first} gives
the main term $\mathcal XG(\ell)$ and a relative error $O(\tau_M)$.
The coefficient bound \eqref{eq:sv-coefficient-bound} gives
\[
 |r_\ell|\ll\mathcal XG(\ell)\tau_M
 +X^{1/2+\eps}M^{2+o(1)}\ell^{1+\eps}.
\]
We have $0\le G(p)<1/p$. To see this, set $x=1/p$ and
$q=Q_p$. Hasse's bound gives
\[
 (1-x)^2\le q\le(1+x)^2,\qquad
 T_p-q=(1-x)((1+x)^2-q)+x^2(1+x)>0.
\]
Moreover $g(p),g(p^2)>0$ and $a(p)b(p)\ge0$, so $A_1(p)\ge1$.
It follows that $0\le G(p)\le\delta(p)<1/p$ when $p\nmid N_0$;
for $p\mid N_0$ this is immediate from $G(p)=0$. In particular,
$B_1\ge1$, and hence $\mathcal X\ge C_1X$.
Thus
\[
 \sum_{\ell\le D}|r_\ell|
 \ll\mathcal X\tau_M\log D
 +X^{1/2+\eps}M^{2+o(1)}D^{2+\eps}.
\]
Indeed, for each sufficiently small fixed $\eps>0$ the exponent of $X$
in the second term is at most
$1/2+2\theta+2\gamma+\eps(1+\gamma)+o(1)$, which is strictly
less than $1$ by $4\theta+4\gamma<1$. Since $\mathcal X\ge C_1X$,
this gives a fixed power saving relative to $\mathcal X$.
The first term is smaller than
$\mathcal X(\log X)^{-A}$ for every fixed $A$.
\end{proof}

\section{The weighted sieve}

\begin{proof}[Proof of Theorem~\ref{thm:main}]
Choose fixed $\rho\in(0,1)$ and $\theta,\gamma>0$ with $c<2\rho\sqrt\theta$. We also require $35/(4(7W+1))<\gamma<1/4-\theta$.
Take the resonator of Section~\ref{sec:resonator}, and put
$D=X^\gamma$, $z=D^{1/5}$,
$P(z)=\prod_{p<z,\,p\nmid N_0}p$ and
$V(z)=\prod_{p\mid P(z)}(1-G(p))$.
The inequality $G(p)<1/p$ gives the uniform dimension-one
condition and $V(z)\gg_E1/\log z$ by Mertens' theorem.

Let $F$ and $f$ be the upper and lower linear sieve functions
\cite[Theorem~11.13]{FI},
and let $C=2e^{\gamma_{\rm EM}}$, where $\gamma_{\rm EM}$ is
Euler's constant. The sieve formulas give
$F(t)=C/t$ for $1<t\le3$ and $f(4)=C\log3/4$.
For $3\le r\le4$, the delay equation gives
\[
 rF(r)=C\left(1+\int_3^r\frac{\log(t-2)}{t-1}\,dt\right)
 \le\frac{Cr}{3},
\]
since $\log(t-2)\le t-3\le(t-1)/3$ on this interval.
It also gives
\[
 5f(5)=4f(4)+\int_3^4F(t)\,dt\ge C\log4,
 \qquad f(5)\ge\frac{C\log4}{5}>f(4).
\]
For $1<v\le4$, put $I(v)=\int_1^v(1-t/v)F(5-t)\,dt/t$.
Splitting at $t=2$ and using the preceding bounds gives
\[
 I(4)\le C\left(\frac8{15}\log2-\frac1{12}
                    -\frac1{20}\log3\right)<\frac67f(4).
\]
Thus $I(v)\le I(4)<(6/7)f(5)$ for $v\le4$.
Fix $\eta=6/7$ and choose $1<v<4$ close enough to $4$ that $\gamma>5/(v(W+1/7))$.
Put $y=z^v$ and define the logarithmic weight
\[
 w(d)=1-\frac1\eta
 \sum_{\substack{z\le p<y\\p\mid d}}
 \left(1-\frac{\log p}{\log y}\right).
\]
For $(d,P(z))=1$, all prime factors of $|d|$ are at least $z$.
If $\omega(|d|)\ge W+1$, then
\[
 \sum_{\substack{z\le p<y\\p\mid d}}
 \left(1-\frac{\log p}{\log y}\right)
 \ge W+1-\frac{\log(2X)}{\log y}>\eta
\]
for large $X$. Thus $w(d)\le0$ on these twists, and
$w(d)\le1$ on all sifted twists.

Write $a_d=L(1/2,E_d)R(d)^2\Phi(|d|/X)\ge0$, and put
\[
 S(z)=\sum_{(d,P(z))=1}a_d,
 \qquad S_p(z)=\sum_{\substack{(d,P(z))=1\\p\mid d}}a_d.
\]
All sums over $d$ below are over $\mathcal D_a^-$.
The lower linear sieve at level $D$ gives
 $S(z)\ge\mathcal XV(z)(f(5)+o(1))$.
For $z\le p<y$, apply the upper linear sieve to the sequence
$\{a_d\mathbf1_{p\mid d}\}$ at level $D/p$.
Its divisor sums are $\mathcal XG(p)G(\ell)+r_{p\ell}$ for
$\ell\mid P(z)$, and its sieve parameter is
$5-\log p/\log z\in[5-v,4]$.
Consequently,
\[
 S_p(z)\le\mathcal XG(p)V(z)
 \left(F\left(5-\frac{\log p}{\log z}\right)+o(1)\right)
 +\sum_{\substack{\ell\le D/p\\\ell\mid P(z)}}|r_{p\ell}|.
\]
The $o(1)$ is uniform in $p$ since $5-v>1$ is fixed.
Each $p\ell$ in this sum has exactly one prime factor at least $z$,
so Lemma~\ref{lem:richert-distribution} gives
\[
 \sum_{z\le p<y}\sum_{\substack{\ell\le D/p\\\ell\mid P(z)}}
 |r_{p\ell}|\le\sum_{q\le D}|r_q|=o(\mathcal XV(z)).
\]
Here $\log\mathcal L=O(\log\log X)$, so $\log P_+=o(\log X)$.
Since $\log z=(\gamma/5)\log X$, we have $P_+<z$ for large $X$.
Also $z>N_0$ eventually. Consequently, for every prime $p\ge z$,
$b(p)=0$, $p\nmid N_0$, and $G(p)=\delta(p)$; the formulas above
and Hasse's bound give the uniform estimate $|G(p)-1/p|\le6/p^2$.
The function $F$ is bounded on the fixed compact interval $[5-v,4]$,
and the logarithmic weight lies in $[0,1]$. Hence replacing $G(p)$
by $1/p$ in the following sum changes it by at most
$O_v(\sum_{n\ge z}n^{-2})=o(1)$. The prime number theorem and partial
summation then give
\[
 \sum_{z\le p<y}\left(1-\frac{\log p}{\log y}\right)
 G(p)F\left(5-\frac{\log p}{\log z}\right)=I(v)+o(1).
\]
The sieve bounds give
\[
 T:=\sum_{(d,P(z))=1}a_dw(d)
 \ge\mathcal XV(z)\left(f(5)-\frac{I(v)}\eta+o(1)\right)
 \gg\frac{XB_1}{\log X}.
\]

Since $T>0$, some term has $w(d)>0$ and $\omega(|d|)\le W$.
Discard the terms with $w(d)\le0$.
For the remaining terms, $w(d)\le1$, so
\[
 T\le\max_{d\in\mathcal F_{a,W}(X)}L(1/2,E_d)
 \sum_{d\in\mathcal D_a^-}R(d)^2\Phi(|d|/X).
\]
The sum on the right is $\ll XB_0$ by
\cite[Proposition~3.1]{HL}, since $\theta<1/4$.
Therefore
\[
 \max_{d\in\mathcal F_{a,W}(X)}L(1/2,E_d)
 \gg\frac{B_1}{B_0\log X}.
\]
By Lemma~\ref{lem:euler-ratio}, $\log(B_1/B_0)=(2\rho\sqrt\theta+o(1))H_X$.
As $\log\log X=o(H_X)$ and $c<2\rho\sqrt\theta$,
the asserted bound follows.
\end{proof}
\section{Large Tate--Shafarevich groups}

For an elliptic curve $A/\mathbb Q$, our normalization gives
$L(1/2,A)=L_{\mathrm{HW}}(1,A)$, where
$L_{\mathrm{HW}}$ is the classical Hasse--Weil $L$-function.
Let $\omega_A$ be a minimal N\'eron differential and put
$\Omega(A)=\int_{A(\mathbb R)}|\omega_A|$ and
$\Tam(A)=\prod_{p<\infty}c_p(A)$.
The integral includes all $c_\infty(A)=\#\pi_0(A(\mathbb R))$
connected components. Write $H_1(A(\mathbb C),\mathbb Z)^\pm$ for
the eigensubgroups of complex conjugation, and define
$\Omega^-(A)=\int_{\gamma_A^-}\omega_A$ using a generator
$\gamma_A^-$ of the minus subgroup. This period is purely imaginary
and determined up to sign.
Put
\begin{equation}\label{eq:analytic-BSD-quantity}
 \mathcal B(A)=\frac{L(1/2,A)\#A(\mathbb Q)_{\mathrm{tors}}^2}
 {\Omega(A)\Tam(A)}.
\end{equation}
The rank-zero Birch--Swinnerton-Dyer formula is
\begin{equation}\label{eq:rank-zero-bsd}
 \#\Sha(A/\mathbb Q)=\mathcal B(A)
 \qquad\text{when }L(1/2,A)\ne0.
\end{equation}

\begin{lemma}\label{lem:arithmetic-factors}
Let $d<0$ be a fundamental discriminant with $(d,2N)=1$.
Then
\begin{equation}\label{eq:negative-period}
 \Omega(E_d)=\frac{C_E^-}{\sqrt{|d|}},
 \qquad C_E^-=c_\infty(E)|\Omega^-(E)|.
\end{equation}
The complementary period satisfies
\[
 |\Omega^-(E_d)|=\frac{\Omega(E)}{c_\infty(E)\sqrt{|d|}},
 \qquad |d|\Omega(E_d)|\Omega^-(E_d)|=\Omega(E)|\Omega^-(E)|.
\]
The Tamagawa product satisfies
\begin{equation}\label{eq:tamagawa-bound}
 1\le\Tam(E_d)\ll_E4^{\omega(|d|)}.
\end{equation}
For any fixed $W\ge1$, we have $\mathcal B(E_d)\gg_{E,W}\sqrt{|d|}\,L(1/2,E_d)$ if $\omega(|d|)\le W$.
\end{lemma}
\begin{proof}
By \cite[Theorem~3.2]{Pal}, $\Omega(E_d)=\widetilde u(d)c_\infty(E_d)|\Omega^-(E)|/\sqrt{|d|}$, where $\widetilde u(d)$ is the minimal-model correction.
Here $d$ is odd, squarefree, and congruent to $1$ modulo $4$.
It is coprime to the minimal discriminant of $E$, so
$\widetilde u(d)=1$ by
\cite[Proposition~2.5 and Corollary~2.6]{Pal}.
Quadratic twisting preserves the sign of the discriminant, so
$c_\infty(E_d)=c_\infty(E)$. This proves
\eqref{eq:negative-period}.

Choose the twisting isomorphism $\phi:E_d\to E$ over $\mathbb C$
so that $\phi^*\omega_E=\sqrt d\,\omega_{E_d}$; see
\cite[Lemma~3.1]{Pal}. Complex conjugation satisfies
$\sigma_E\circ\phi=[-1]\circ\phi\circ\sigma_{E_d}$.
Thus $\phi_*$ exchanges the two integral homology eigensubgroups
and maps $\gamma_{E_d}^-$ to a generator of
$H_1(E(\mathbb C),\mathbb Z)^+$.
The absolute period of this generator is
$\Omega(E)/c_\infty(E)$; see \cite[Lemma~5.1]{Pal}. Hence
\[
 \sqrt{|d|}\,|\Omega^-(E_d)|
 =\left|\int_{\phi_*\gamma_{E_d}^-}\omega_E\right|
 =\frac{\Omega(E)}{c_\infty(E)}.
\]
The complementary formula and the product identity follow.

We have $c_p(E_d)=1$ for $p\nmid dN$.
If $p\mid d$, then $E$ has good reduction at $p$ and $E_d$
has type $I_0^*$, with
$c_p(E_d)=\#E(\mathbb F_p)[2]\in\{1,2,4\}$; see
\cite[Theorem~4.1 and Table~4]{Local}.
For $p\mid2N$, the local twist depends only on the squareclass of
$d$ in $\mathbb Q_p^\times$. There are finitely many such classes,
so their Tamagawa numbers are bounded in terms of $E$.
This gives \eqref{eq:tamagawa-bound}.
Since $\#E_d(\mathbb Q)_{\mathrm{tors}}\ge1$, the last bound follows
from \eqref{eq:negative-period}, \eqref{eq:tamagawa-bound} and
\eqref{eq:analytic-BSD-quantity}.
\end{proof}

\begin{proof}[Proof of Corollary~\ref{cor:Sha-main}]
Take $c'$ such that $c<c'<\sqrt{(7W-34)/(7W+1)}$.
Theorem~\ref{thm:main} gives
$d\in\mathcal F_{a,W}(X)$ with $L(1/2,E_d)\ge\exp(c'H_X)$.
Lemma~\ref{lem:arithmetic-factors} gives
$\mathcal B(E_d)\gg_{E,W}\sqrt X\exp(c'H_X)$.
Since $c'>c$, we have $\mathcal B(E_d)\ge\sqrt X\exp(cH_X)$ for all sufficiently large $X$.

The nonzero central value implies that $E_d(\mathbb Q)$ and
$\Sha(E_d/\mathbb Q)$ are finite by the
Gross--Zagier--Kolyvagin theorem \cite{GZ,Kol}.
Under \eqref{eq:rank-zero-bsd}, this proves \eqref{eq:Sha-large}.
If $E$ has complex multiplication, then so does $E_d$.
By \cite[Corollary~2 and the paragraph following its proof, p.~362]{BF},
the full rank-zero BSD formula holds for every CM elliptic curve over
$\mathbb Q$ with nonzero central value. This gives the unconditional
assertion for CM curves.
\end{proof}

\begin{remark}[The Tamagawa denominator]
For $d\in\mathcal D_a^-$, the local squareclass of $d$ is fixed at every
$p\mid2N$. The number $B_{E,a}=\prod_{p\mid2N}c_p(E_d)$ is therefore
constant, and $\Tam(E_d)=B_{E,a}\prod_{p\mid d}\#E(\mathbb F_p)[2]$.
Frobenius acting on $E[2]$ determines each factor; see also
\cite[\S5]{RS}. We have
\[
 \max_{\substack{d\in\mathcal D_a^-\\X\le|d|\le2X}}
       \log\Tam(E_d)
 = (\log4+o(1))\frac{\log X}{\log\log X}.
\]
For the upper bound, use $\Tam(E_d)\le B_{E,a}4^{\omega(|d|)}$.
If $r=\omega(n)$, then $n\ge(r+1)!$. Stirling's formula gives
$\omega(n)\le(1+o(1))\log X/\log\log X$ uniformly for $n\le2X$.

For the lower bound, let $K=\mathbb Q(E[2])$ and
$\delta=[K:\mathbb Q]^{-1}$. Let $\mathcal P$ be the primes $p\nmid2N$
that split completely in $K$, and put
$\pi_{\mathcal P}(y)=\#\{p\le y:p\in\mathcal P\}$.
Each such prime satisfies $\#E(\mathbb F_p)[2]=4$. Chebotarev gives
\[
 \pi_{\mathcal P}(y)\sim\frac{\delta y}{\log y},
 \qquad \sum_{\substack{p\le y\\p\in\mathcal P}}\log p\sim\delta y;
\]
see \cite[Lemma~2.4]{LM}. Fix $0<\varepsilon<1$, and put
$y=(1-\varepsilon)\log X/\delta$ and
$m=\prod_{p\le y,\,p\in\mathcal P}p$.
Then $\log m=(1-\varepsilon+o(1))\log X$ and
$Y=X/m=X^{\varepsilon+o(1)}$.
The prime number theorem in fixed arithmetic progressions gives
a prime $q\in[Y,2Y]$ with $q\equiv-am^{-1}\pmod{N_0}$.
For large $X$, we have $q>y$ and $q\nmid m$.
The integer $d=-qm$ is squarefree, belongs to $\mathcal D_a^-$, and
satisfies $X\le|d|\le2X$. Therefore
\[
 \log\Tam(E_d)\ge\pi_{\mathcal P}(y)\log4
 =\bigl((1-\varepsilon)\log4+o(1)\bigr)
        \frac{\log X}{\log\log X}.
\]
Let $\varepsilon\to0$ to obtain the lower bound.

If $E(\mathbb Q)[2]\ne0$, the local product gives
$\Tam(E_d)\ge B_{E,a}2^{\omega(|d|)}$.
If all $2$-torsion points are rational, then
$\Tam(E_d)=B_{E,a}4^{\omega(|d|)}$.
For these curves, a uniform bound for $\Tam(E_d)$ also bounds
$\omega(|d|)$. The almost prime condition controls this denominator
for every fixed elliptic curve.
More generally, the condition $\log\Tam(E_d)=o(H_X)$ on the selected
twists preserves the central-value gain.
\end{remark}

\section*{Acknowledgements}
S.H. would like to thank Chung Pang Mok for the constant encouragement.
S.H. is supported by the National Natural Science Foundation of China (No. 12601008) and the China Postdoctoral Science Foundation (No. 2026M793348).

\section*{AI Use Disclosure}
The author introduced sieve weights into the argument, using calculations based on earlier work, and first reduced the number of prime factors to six. ChatGPT introduced Richert's logarithmic weighted sieve to reduce the bound to five and prepared the draft.

\end{document}